\documentclass{article}

\usepackage{arxiv}
\usepackage[utf8]{inputenc} 
\usepackage[T1]{fontenc}    
\usepackage{hyperref}       
\usepackage{url}            
\usepackage{booktabs} 
\usepackage{amsfonts} 
\usepackage{nicefrac} 
\usepackage{microtype}      
\usepackage{graphicx}
\usepackage{doi}
\usepackage{subcaption}
\usepackage{amsmath,amsfonts,amssymb,mathtools}
\usepackage{multirow}
\usepackage{amsthm}
\usepackage{xcolor}

\usepackage{tabularx}
\usepackage{array}
\usepackage{threeparttable}
\usepackage{enumitem}
\usepackage{hyperref}
\hypersetup{hidelinks}
\usepackage[nameinlink,noabbrev]{cleveref}
\usepackage{placeins}
\usepackage{float}
\usepackage{flafter}
\usepackage{longtable}
\usepackage{adjustbox}
\usepackage{needspace}

\title{Exterior-Row Replacement for Grid-Based Multinode Shepard Collocation on Implicit Planar Domains}

\date{}

\author{ \href{https://orcid.org/0000-0001-7608-5635}{\includegraphics[scale=0.06]{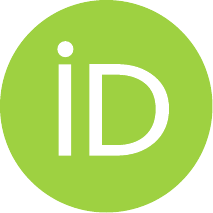}\hspace{1mm}Anouar El~Harrak}\thanks{Corresponding author.} \\
	SMMADCSN Laboratory, FPL \\
    Abdelmalek Essaadi University\\
	Tetouan, Morocco \\
	\texttt{a.elharrak@uae.ac.ma} \\
    \And
    \href{}{\includegraphics[scale=0.06]{orcid.pdf}\hspace{1mm}Amal Bergam} \\
	SMMADCSN Laboratory, FPL \\
    Abdelmalek Essaadi University\\
	Tetouan, Morocco \\
	\texttt{abergam@uae.ac.ma} \\
    \And
	\href{https://orcid.org/0000-0002-5988-650X}{\includegraphics[scale=0.06]{orcid.pdf}\hspace{1mm}Hatim Tayeq} \\
	SMMADCSN Laboratory, FPL \\
    Abdelmalek Essaadi University\\
	Tetouan, Morocco \\
	\texttt{h.tayeq@uae.ac.ma} \\
}

\renewcommand{\shorttitle}{}
\newtheorem{theorem}{Theorem}[section]
\newtheorem{proposition}[theorem]{Proposition}
\newtheorem{lemma}[theorem]{Lemma}

\theoremstyle{definition}

\newcommand{\R}{\mathbb{R}}

\newcommand{\I}{\mathcal{I}}
\newcommand{\Sset}{\mathcal{S}}

\newcommand{\rank}{\operatorname{rank}}
\newcommand{\dist}{\operatorname{dist}}
\newcommand{\RMS}{\mathrm{RMS}}

\newcommand{\epsmach}{\varepsilon_{\mathrm{mach}}}
\newcommand{\sigmin}{\sigma_{\min}}
\newcommand{\kappaTwo}{\kappa_2}
\newcommand{\qandq}{\quad\text{ and }\quad}

\begin{document}
\maketitle

\begin{abstract}
Grid-based multinode Shepard collocation combines tensor-product Lagrange interpolants on overlapping Cartesian subgrids. On an implicitly defined domain, retained subgrids may cross the boundary and contain exterior source nodes. Removing or moving these nodes changes the local interpolants. We introduce exterior-row replacement (ERR) to preserve the source nodes and coefficients of every retained subgrid. We assemble PDE rows at strictly interior sources and supply the remaining equations through boundary evaluations. We consider two square completion methods based on projection and redistribution, an oversampled weighted least-squares one (ERR-OS), and a completion method that selects boundary rows through the interior PDE nullspace (ERR-QR). Full column rank is characterized by the boundary action on this nullspace. We also give completion conditions for prescribed fitted rows and an error estimate in terms of the residual. Global-full (GF) support uses all retained subgrids, whereas nearest-active (NA) support uses only those containing the active source nearest to the evaluation point. NA keeps off-source boundary rows sparse. The numerical study covers a disk, an annulus, a smooth star-shaped domain, a square hole with fitted and implicit boundary components, and a fitted L-shaped domain. With degree-$4$ local tensor polynomials, ERR-OS and ERR-QR recover the expected fifth-order accuracy in the reported refinement study. Further experiments examine the effects of boundary oversampling and trimming on algebraic admissibility and boundary observability.
\end{abstract}


\Needspace{6\baselineskip}
\section{Introduction}
\label{sec:introduction}

Cartesian source grids provide a simple data structure for constructing local approximations and performing numerical differentiation \cite{JOHANSEN199860}. The boundary of the domain typically does not coincide with the computational grid \cite{burman2025cutfem, tominec2021unfitted}. A numerical method must therefore distinguish the nodes defining the source representation from the points at which the PDE and boundary conditions are collocated. Immersed finite elements \cite{deprenter2023immersed, https://doi.org/10.1002/nme.4823, HANSBO20025537}, embedded and shifted-boundary methods \cite{main2018shifted}, ghost-point schemes \cite{coco2013ghost}, and unfitted meshfree approaches address this geometric mismatch using different constructions \cite{nielsen2025unfitted,gibou2002poisson,bayona2017rbffd,larsson2017pum,tominec2021lsrbffd}.
   
Multinode Shepard collocation has been used for Poisson and other elliptic Dirichlet problems on domains of different geometries \cite{dellaccio2021poisson,dellaccio2024domains}. Enriched multinode Shepard spaces have also been considered for elliptic problems with singular boundary data \cite{dellaccio2024enriched}. Recent studies have developed multinode Shepard functions based on tensor-product interpolation over uniform rectangular grids \cite{barrera2026tensor}. The separation between source and collocation sets is also established in unfitted collocation \cite{doi:10.1137/21M142126X}. The least-squares RBF-FD method proposed in \cite{tominec2021unfitted} uses an unfitted node structure and additional collocation points. The formulation presented in \cite{nielsen2025unfitted} employs interior nodes and enforces boundary constraints via Lagrange multipliers. Boundary treatment and row replacement are standard issues in spectral collocation. In this context, rectangular spectral collocation similarly addresses the balance between differential and boundary equations \cite{driscoll2016rectangular}.
   
Grid-based multinode Shepard collocation (GBMSC) combines Lagrange tensor-product interpolants defined on overlapping Cartesian subgrids \cite{elharrak2026gbmsc}. For an implicit domain, a background grid covers the domain and a subgrid is retained when its fraction of strictly interior nodes exceeds a prescribed threshold. A retained subgrid may extend outside the domain so that certain nodes of the subgrid lie outside the domain. These exterior nodes are retained as part of the local interpolation subgrid. The elimination of these exterior nodes modifies and alters the Cartesian structure of those subgrids near the boundary. In contrast, imposing the partial differential equation (PDE) at an exterior source node is inadmissible within the framework of the boundary-value problem. The present work concerns the assembly of boundary completion methods that preserve each retained complete tensor grid while keeping PDE collocation points only in the interior of the domain. We address this assembly by preserving all active sources of the retained subgrids and their associated coefficients. The PDE is collocated only at active sources classified as strictly interior. The remaining equations are supplied by boundary evaluations, whose collocation points need not be source nodes. This procedure is called exterior-row replacement (ERR).
   
We consider four completion methods with different objectives. ERR-P maps one projected boundary point to each exterior source point and then uses midpoint reprojection to resolve duplicate points. This approach produces a square system provided that the geometric procedure yields the necessary number of boundary points. The second completion method, ERR-U, uses a uniform distribution of boundary points, whereas the completion method, ERR-OS$(\rho)$, employs an oversampling of points on the boundary with a ratio $\rho>1$ and solves a weighted least-squares problem \cite{cheung2018h}. In the last completion method, ERR-QR, we select the boundary rows based on their action on the nullspace of the interior PDE block. First, we study the rank conditions and derive a stability estimate. In fact, the rank identity follows from the rank-nullity, and the selection procedure uses the classical column-pivoted QR \cite{businger1965householder}. Next, we apply these results to the ERR completion methods combined with fitted rows. The polynomial reproduction introduced in the study ensures that sampled harmonic polynomials belong to the interior PDE nullspace, which gives a lower bound on the dimension of this nullspace. The stability estimate provides a bound on the smallest singular value of the full system by combining the contribution of the interior PDE operator, the action of the boundary operator restricted to the nullspace, and the coupling between these components.
   
The numerical experiments compare the four completion methods under refinement on a disk, an annulus, and a smooth star-shaped domain. We examine the oversampling ratio, trimming, and the support policy. In the comparison experiment between ERR-OS and ERR-QR, we employ the same boundary points to quantify the impact of retaining or compressing redundant boundary information. A hybrid square hole and a fitted L-shaped domain examine the use of fitted and implicit boundary rows in the same assembly framework. 
    
The remainder of the paper is structured as follows. Section~\ref{sec:problem} introduces the governing PDE, the associated implicit geometric representation, and the Cartesian GBMSC components required for the present study. Section~\ref{sec:err} specifies the retained source space, the common boundary-row framework, and the four completion methods considered. Section~\ref{sec:analysis} presents the geometric and algebraic analyses. Section~\ref{sec:experiments} reports the numerical experiments.

\Needspace{6\baselineskip}
\section{Grid-based multinode Shepard collocation on implicit domains}
\label{sec:problem}

Let $\Omega\subset\mathbb R^2$ be a bounded open domain with boundary $\Gamma$. We study the Dirichlet boundary value problem
\begin{equation}
\begin{cases}
-\Delta u = f, & \text{in } \Omega,\\[2mm]
u = g, & \text{on } \Gamma.
\end{cases}
\label{eq:poisson}
\end{equation}

Here, $u$ is the unknown scalar function, and $f$ and $g$ are the prescribed source and Dirichlet data. The geometry is represented by a sign function $\phi$ satisfying
\[
\phi<0\ \hbox{in }\Omega,\qquad
\phi=0\ \hbox{on }\Gamma,\qquad
\phi>0\ \hbox{outside }\overline\Omega.
\]
The algebraic construction uses the source classification and the supplied boundary collocation points. Regularity of $\phi$ is required locally where the retraction is used. If $\phi$ is the exact signed distance in a regular tubular neighborhood, then $|\nabla\phi|=1$ in that neighborhood and, for exterior points, $\phi(x)=\dist(x,\Omega)$.
   
Given a starting point $x^{(0)}=x$, the numerical retraction uses the Newton-type iteration
\[
x^{(k+1)}=x^{(k)}-
\phi(x^{(k)})\frac{\widetilde\nabla_h\phi(x^{(k)})}
{\|\widetilde\nabla_h\phi(x^{(k)})\|_2^2},
\]
where $\widetilde\nabla_h\phi$ denotes the finite-difference gradient used in the implementation. When the iteration returns a point that meets the prescribed zero-level tolerance, that point is denoted by $R_\Gamma(x)$. This point is an approximation to a boundary point.
     
Let $G_h$ be a uniform Cartesian grid with spacing $h>0$ on a rectangular box containing $\Omega$. The construction of the GBMSC method is described in \cite{elharrak2026gbmsc}. We use isotropic tensor subgrids with $p+1$ consecutive coordinates in each direction, where $p\ge2$. Thus, each $X_s$ has $(p+1)^2$ nodes and diameter at most $\sqrt2ph$. The completion methods could also be applied to other tensor layouts, but the geometric constants and trimming fractions must then be adjusted. We define the fraction of the strictly interior nodes of a subgrid $X_s$ by
\[
\vartheta_s=\frac{\#\{z\in X_s:\phi(z)<0\}}{(p+1)^2}.
\]
For $\vartheta_{\min}$ the minimum fraction of strictly interior nodes, $0\le\vartheta_{\min}<1$, the retained family of subgrids $\Sset_h^{\rm ret}$ and the active source set $X_h^{\rm act}$ are
\begin{equation}
\Sset_h^{\rm ret}=\{s:\vartheta_s>\vartheta_{\min}\}
\qandq X_h^{\rm act}=\bigcup_{s\in\Sset_h^{\rm ret}}X_s.
\label{eq:retained_family}
\end{equation}
We assume that $\Sset_h^{\rm ret}$ is not empty. Every retained subgrid contains at least one strictly interior source. A retained subgrid that crosses the boundary continues to constitute a complete tensor-product subgrid. In particular, any exterior node contained within a retained crossing subgrid remains an active source coefficient.

Given nodal data $v_z$ of a scalar function $v$ at the node $z$, the corresponding tensor-product Lagrange interpolant is defined as
\begin{equation*}
L_s [v](x)=\sum_{z\in X_s}v_z\,\ell_{s,z}(x).
\end{equation*}

$\ell_{s,z}$ is the tensor-product Lagrange basis function associated
with the node $z$ on the Cartesian subgrid $X_s$. These functions satisfy the cardinal property
\[
\ell_{s,z}(z')=\delta_{zz'},
\qquad z,z'\in X_s.
\]

Away from active source nodes, the unnormalized and normalized multinode Shepard weights are, respectively,
\begin{equation*}
q_s(x)=\prod_{z\in X_s}\|x-z\|^{-\mu}
\qandq
w_s(x)=\frac{q_s(x)}{\sum_{t\in\Sset_h(x)}q_t(x)}.
\end{equation*}
  
We assume $\mu>2$ to ensure that the solution is twice continuously differentiable at interior source points, as required for PDE collocation. At the source nodes, the singular expressions are interpreted through their limiting form. The contributor set $\mathcal S_h(x)$ is defined by one of the following two support policies:

\begin{itemize}

\item \textbf{Global-full policy (GF)}: GF uses every retained contributing subgrid. Let
    $$
    \mathcal S_h^{\rm GF}(x)=\mathcal S_h^{\rm ret}.
    $$

\item \textbf{Nearest-active policy (NA)}: NA uses subgrids associated with the nearest active source. Let
$$
x^*\in
\operatorname*{arg\,min}_{z\in X_h^{\rm act}}\|x-z\|_2.
$$
The contributor set of NA policy is defined as the set of retained subgrids incident to the nearest active source. Let
$$
\mathcal S_h^{\rm NA}(x)
=
\{\,s\in\mathcal S_h^{\rm ret}:x^*\in X_s\,\}.
$$

\end{itemize} 

The set $\mathcal S_h(x)$ denotes the set of GF or NA contributors. For NA, we randomly select an active source when multiple sources are located at comparable distances. GF has a $C^2$ nodal extension for $\mu>2$. NA is smooth within each nearest-source region, but it is not necessary to be continuous across the interfaces.
   
For a set of nonempty contributors $\Sset_h(x)$, the global approximation is
\begin{equation*}
I_h[v](x)=\sum_{s\in\Sset_h(x)}w_s(x)L_s[v](x).
\end{equation*}
We approximate the solution of the Poisson problem \eqref{eq:poisson} by
\begin{equation*}
u_h(x)=I_h[c_h](x),
\end{equation*}
where $c_h=\{c_z\}_{z\in X_h^{\rm act}}$ is the vector of unknown source coefficients. For a support policy and at off-source points, define the cardinal functions and the Dirichlet evaluation row by
\begin{equation}
\psi_z(x)=\sum_{\substack{s\in\mathcal S_h(x)\\z\in X_s}}w_s(x)\ell_{s,z}(x)
\qandq
\mathcal B_h(y)=(\psi_z(y))_{z\in X_h^{\rm act}}.
\label{eq:cardinal_rows}
\end{equation}
Then $I_h[c](x)=\sum_z\psi_z(x)c_z$ and $\mathcal B_h(y)c=I_h[c](y)$. The same formulas use the nodal limits at active sources. GF and NA have the same values and derivatives up to the second order at each active source when $\mu>2$. Thus, the nodal PDE rows coincide exactly at source nodes, whereas their evaluations at off-source points generally differ.

\section{Exterior-row replacement and boundary completion methods}
\label{sec:err}

Let $\tau_{\rm in}\ge0$ be the sign tolerance on the grid $G_h$. We partition the set of active indices into
\begin{equation}
\I_-=\{j:\phi(x_j)<-\tau_{\rm in}\},
\label{eq:Iminus}
\end{equation}
\begin{equation}
\I_0=\{j:|\phi(x_j)|\le \tau_{\rm in}\},
\label{eq:Izero}
\end{equation}
\begin{equation}
\I_+=\{j:\phi(x_j)>\tau_{\rm in}\},
\qandq
\I_{\rm rep}=\I_0\cup\I_+.
\label{eq:Iplus}
\end{equation}
Let
\begin{equation*}
n_h=|X_h^{\rm act}|,
\qquad
n_-=|\I_-|
\qandq
n_{\rm rep}=|\I_{\rm rep}|=n_h-n_-.
\end{equation*}
The PDE is assembled exclusively at the source locations indexed by $\I_-$. The corresponding interior PDE block and the associated right-hand-side vector are given by
\begin{equation*}
P_h\in\R^{n_-\times n_h} \qandq f_h
=
\bigl(f(x_i)\bigr)_{i\in\mathcal I_-}
\in\mathbb R^{n_-}.
\end{equation*}
For boundary collocation points $\{y_j\}_{j=1}^{m_\Gamma}$, the rows corresponding to the boundary evaluations and the associated right-hand-side vector are given by
\begin{equation*}
B_h\in\R^{m_\Gamma\times n_h} \qandq g_h
=
\bigl(g(y_j)\bigr)_{j=1}^{m_\Gamma}
\in\mathbb R^{m_\Gamma}.
\end{equation*}
If a retracted collocation point is not exactly on $\Gamma$, the evaluation of the boundary data $g$ at that point requires extending $g$ to a neighborhood of $\Gamma$.
   
The complete system can be expressed as a block-matrix as follows.
\begin{equation*}
A_hc_h
=
\begin{bmatrix}P_h\\B_h\end{bmatrix}c_h
=
\begin{bmatrix}f_h\\g_h\end{bmatrix}=b_h.
\end{equation*}

For an implicit boundary component, the numerical retraction procedure provides the location of the boundary collocation point. For a fitted boundary component, the fitted sources are used directly as collocation points. Duplicate or near-duplicate points are treated by the repair reprojection procedure described below.

\subsection{Boundary completion methods for unfitted boundary components}
We consider four completion formulations, namely ERR-P, ERR-U, ERR-OS, and ERR-QR. Each addresses a different methodological objective. The ERR-U, ERR-OS, and ERR-QR constructions share a denser finite cloud of candidate points
\begin{equation}
Y_h^{\rm cand}=\{y_i^{\rm cand}\}_{i=1}^{M_h}\subset\Gamma,
\qquad
M_h\ge n_{\rm rep}.
\label{eq:candidate_cloud}
\end{equation}
The corresponding candidate boundary matrix is defined as follows:
\begin{equation*}
C_h
=
\begin{bmatrix}
\mathcal B_h(y_1^{\rm cand})\\
\vdots\\
\mathcal B_h(y_{M_h}^{\rm cand})
\end{bmatrix}
\in\R^{M_h\times n_h}.
\end{equation*}
Here $\mathcal B_h(y)$ is the boundary row defined in \eqref{eq:cardinal_rows}.
   
\subsubsection*{ERR-P: projection-based replacement}
\label{sec:errp}

In the ERR-P procedure, one boundary point is associated with each exterior source point $x_j$, where $j\in\I_{+}$. Once the zero-level residual satisfies the specified tolerance criterion, the returned boundary point is defined as
\[
y_j=R_\Gamma(x_j).
\]
   
Different sources may be retracted to the same (or numerically indistinguishable) point. If $y_j$ is too close to an existing boundary point, we choose a neighboring distinct point $y_\ell$ and form the midpoint
\[
m_{j\ell}=\tfrac12(y_j+y_\ell).
\]
The midpoint is reprojected, $\widetilde y_j=R_\Gamma(m_{j\ell})$, repeated as needed until a distinct boundary point is found. The new point replaces $y_j$ as the boundary collocation point for the source node $x_j$. If executed successfully, ERR-P retains one boundary equation per replaced exterior source point,
\begin{equation*}
m_\Gamma=n_{\rm rep}.
\end{equation*}
With the $n_-$ interior PDE equations, we get the $n_h$ equations for the $n_h$ active Cartesian coefficients. This projection with repairing constructs only the boundary collocation point to complete the square system, but the exterior source point and the retained subgrids are unchanged. In the ERR-P completion method, we construct boundary rows without using rank information from the interior PDE block.

\subsubsection*{ERR-U: redistributed quasi-uniform square completion method}
\label{sec:erru}

ERR-U likewise enforces
\begin{equation*}
m_\Gamma=n_{\rm rep},
\end{equation*}
but selects the associated rows from the denser candidate cloud \eqref{eq:candidate_cloud} according to the prescribed quasi-uniform redistribution rule. The goal is geometric; we select a subset of $Y_h^{\rm cand}$ that covers each boundary component more uniformly while satisfying the square system.

ERR-U may be regarded as the zero-oversampling limiting case of the quasi-uniform sampling strategy, which is also adopted by ERR-OS as described below. This distinction is motivated by experimental design considerations. Specifically, ERR-U isolates the effect of geometric resampling in the absence of additional boundary equations and tests whether such redistribution alone is sufficient.

\subsubsection*{ERR-OS$(\rho)$: quasi-uniform oversampled completion method}
\label{sec:erros}

Let $\rho>1$ denote the boundary oversampling ratio. The ERR-OS$(\rho)$ procedure selects
\begin{equation*}
m_\Gamma(\rho)=\left\lceil\rho n_{\rm rep}\right\rceil
\end{equation*}
quasi-uniformly distributed boundary points from the candidate set $Y_h^{\rm cand}$, subject to the constraint $m_\Gamma(\rho)\le M_h$. We solve the overdetermined system by weighted least squares. For each non-null row, both the row and its corresponding right-hand-side entry are scaled by the Euclidean norm of that row. More precisely, we set two diagonal matrices such that
\[
(D_P)_{ii}=\|(P_h)_{i,:}\|_2^{-1}\qandq
(D_B)_{jj}=\|(B_h)_{j,:}\|_2^{-1}.
\]
For $n_->0$ and $m_\Gamma>0$, ERR-OS minimizes
\begin{equation}
\frac1{n_-}\|D_P(P_hc_h-f_h)\|_2^2+
\frac1{m_\Gamma}\|D_B(B_hc_h-g_h)\|_2^2.
\label{eq:weighted_ls}
\end{equation} 
The objective function is computed as the mean of the squared normalized residuals within each block. Repeating every row and its corresponding right-hand-side entry the same number of times in one block does not alter the block’s mean value. Since $\rho$ is not uniquely defined, the study tests several values, using $\rho = 1.5$ as a reference for comparison.

\subsubsection*{ERR-QR: selection of boundary rows}
\label{sec:errqr}
The ERR-QR boundary completion uses the candidate matrix $C_h$. ERR-QR selects rows from the matrix $C_h$ based on their action on $\ker P_h$. Let
\begin{equation*}
k_h=n_h-\rank(P_h),
\end{equation*}
and let $Z_h\in\R^{n_h\times k_h}$ have orthonormal columns that span $\ker P_h$. The matrix
\begin{equation*}
W_h=C_hZ_h\in\R^{M_h\times k_h}
\end{equation*}
represents the candidate boundary evaluations on this nullspace.

If the condition $\rank(W_h)=k_h$ is satisfied, $C_h$ is injective onto $\ker P_h$. We perform a column-pivoted QR factorization of $W_h^T$ and select the candidate rows from $C_h$ corresponding to the first $k_h$ pivots. These rows define a nonsingular restriction of $\ker P_h$ if $M_h\ge k_h$; see Proposition~\ref{prop:qr_completion}. If $P_h$ has full row rank, then $k_h = n_{\rm rep}$ and the stacked system is square.
   
\section{Geometry and boundary observability}
\label{sec:analysis}

We begin this section with the analysis of the retained Cartesian source space. We then consider the action of the boundary rows on the nullspace of the interior PDE block.

\subsection{Retained sources and polynomial consistency}
\label{sec:retained_analysis}

\begin{lemma}
\label{lem:localization}
\label{prop:count}
Let every retained subgrid $X_s$ contain at least one strictly interior source and consist of $p+1$ consecutive nodes in each direction of a Cartesian grid of spacing $h$. Then every strictly exterior active source $x_j$ satisfies
\begin{equation}
0<\dist(x_j,\Omega)\le\sqrt{2}\,ph.
\label{eq:localization}
\end{equation}
If $\phi$ is the exact signed-distance function, the distance in \eqref{eq:localization} equals $\phi(x_j)$. If, in addition, $\Gamma$ has finite length and a regular tubular neighborhood, then, for fixed $p$ and sufficiently small $h$,
\begin{equation}
|\I_+|=O(h^{-1}).
\label{eq:count}
\end{equation}
\end{lemma}

\begin{proof}
We choose a retained subgrid $X_s$ containing $x_j$ and a strictly interior node $z$ in that subgrid. The diameter of $X_s$ is at most $\sqrt{2}ph$. Consequently,
\[
\dist(x_j,\Omega)\le\|x_j-z\|_2\le\sqrt{2}ph.
\]
If $\phi$ is the exact signed-distance function, we have $$\phi(x_j)=\dist(x_j,\Omega)\le\sqrt{2}\,ph.$$

We choose \(h\) sufficiently small such that the grid cells centered at the exterior sources are contained within the regular tubular neighborhood of \(\Gamma\). We associate with each exterior grid source its square Cartesian cell of area $h^2$. These cells have disjoint interiors. Every such cell lies in the two-sided boundary neighborhood of width $(\sqrt{2}p+\sqrt{2}/2)h$. The regular tubular neighborhood and finite boundary length imply that this neighborhood has area at most $C_\Gamma h$. Hence the number of exterior active sources times $h^2$ is at most $C_\Gamma h$. The set $\I_+$ is a subset of these sources, which proves \eqref{eq:count}.
\end{proof}

The estimate concerns exterior sources. It does not automatically give the same order for the full replacement budget $n_{\rm rep}$ when the interior sign tolerance is fixed. For a signed-distance function and a tolerance $\tau_{\rm in}$ within the regular tubular neighborhood, the same cell argument gives
\[
n_{\rm rep}=O(h^{-1}+\tau_{\rm in}h^{-2}).
\]
Indeed, the tolerance level adds an interior margin of width $\tau_{\rm in}$. Thus, $n_{\rm rep}=O(h^{-1})$ follows if $\tau_{\rm in}=O(h)$. The experiments use the fixed tolerance in Table~\ref{tab:numerical_parameters}.

Under the strict criterion \eqref{eq:retained_family},
\begin{equation}
\vartheta_{\min}\ge 1-(p+1)^{-2}
\label{eq:high_trim}
\end{equation}
implies that every retained subgrid contains only strictly interior sources. Any subgrid containing exactly one source that is not in the interior has an interior fraction bounded above by $1-(p+1)^{-2}$; consequently, such a subgrid cannot be retained. Consequently, $\I_+=\varnothing$. Source nodes with indices in $\I_0$ require separate treatment because, even when a source node is strictly interior, it may still fall within the prescribed sign-tolerance interval.

\begin{proposition}
\label{prop:polynomial_reproduction}
Let $p\ge2$ be an integer and $\mu>2$. Assume that the retained family $\Sset_h^{\rm ret}$ is nonempty and that each retained subgrid $X_s$ is a complete tensor grid with $p+1$ nodes in each direction. Define
\[
Q_p=\operatorname{span}\{x_1^i x_2^j:0\le i,j\le p\},
\qquad
E_h:Q_p\longrightarrow\mathbb R^{n_h},
\qandq
E_hq=(q(z))_{z\in X_h^{\rm act}}.
\]
Then $E_h$ is injective and
\[
I_h[E_hq]=q
\qquad\text{for every }q\in Q_p.
\]

Let $P_h$ contain exact nodal evaluations of $-\Delta I_h[c]$, and let $B_h$ contain exact Dirichlet evaluations at points on $\Gamma$. Suppose that the exact solution belongs to $Q_p$, the data are evaluated exactly, and
\[
A_h=\begin{bmatrix}P_h\\B_h\end{bmatrix}
\]
has full column rank. Then the ERR solution reproduces the exact solution, both for a square system and for weighted least squares with positive weights.
\end{proposition}

\begin{proof}
Suppose that $E_hq=0$ for some $q\in Q_p$. Choose a retained subgrid
\[
X_s=\{\xi_0,\ldots,\xi_p\}\times\{\eta_0,\ldots,\eta_p\}.
\]
All its nodes belong to $X_h^{\rm act}$, so
\[
q(\xi_i,\eta_j)=0,
\qquad 0\le i,j\le p.
\]
For fixed $j$, the polynomial $q(\,\cdot\,,\eta_j)$ has degree at most $p$ and vanishes at the $p+1$ distinct points $\xi_i$. Hence $q(\xi_i,\eta_j)=0$ for every $\xi_i$. Now fix $\xi_i$. The polynomial $q(\xi_i,\,\cdot\,)$ vanishes at the $p+1$ distinct points $\eta_j$ and also has degree at most $p$. It is therefore zero. This proves $q=0$ and the injectivity of $E_h$.

Each local tensor Lagrange interpolant reproduces $q$. Since the weights sum to one,
\[
I_h[E_hq](x)
=\sum_{s\in\mathcal S_h(x)}w_s(x)q(x)
=q(x).
\]
The identity holds for off-source nodes and extends to the source nodes by the nodal limits. For NA, it also holds at a nearest-source node, since every selected local interpolant equals $q$.

Assuming that the exact solution satisfies $q \in Q_p$, the corresponding coefficient vector $c^\star = E_h q$ satisfies
\[
P_hc^\star=f_h
\qandq
B_hc^\star=g_h.
\]
The full column rank makes it the unique solution of a consistent square system and the unique minimizer of weighted least squares. Its reconstruction is $q$.
\end{proof}

\begin{proposition}
\label{prop:polynomial_kernel}
Assuming the hypotheses of Proposition~\ref{prop:polynomial_reproduction}. We define
\[
\mathcal H_p=\{q\in Q_p:\deg q\le p,\ \Delta q=0\}
\qandq
k_h=\dim\ker P_h,
\]
where $\deg q$ denotes the total degree of $q$. Then
\begin{equation}
k_h\ge\dim\mathcal H_p=2p+1.
\label{eq:harmonic_nullity}
\end{equation}
If the $m_\Gamma$ rows of the boundary matrix $B_h$ make
\[
\begin{bmatrix}P_h\\B_h\end{bmatrix}
\]
full column rank, then $m_\Gamma\ge k_h$. The full-row
rank of $P_h$ implies
\[
n_{\rm rep}=n_h-n_-\ge2p+1.
\]
\end{proposition}

\begin{proof}
For $j\ge1$, write
\[
q_j=\sum_{\ell=0}^j a_\ell x_1^{j-\ell}x_2^\ell.
\]
The condition $\Delta q_j=0$ gives
\[
(j-\ell)(j-\ell-1)a_\ell
+(\ell+2)(\ell+1)a_{\ell+2}=0,
\qquad 0\le\ell\le j-2.
\]
Therefore, the coefficients $a_0$ and $a_1$ uniquely determine all other coefficients. Furthermore, the polynomials
\[
\Re(x_1+\mathrm{i}x_2)^j
\qandq
\Im(x_1+\mathrm{i}x_2)^j
\]
are harmonic and linearly independent. Hence, for every degree $j\ge 1$, the associated space of homogeneous harmonic functions is two-dimensional. The property of harmonicity is satisfied independently by each homogeneous component and also by the constants. So, we get
\[
\dim\mathcal H_p=1+2p.
\]

For $q\in\mathcal H_p\subset Q_p$, the polynomial reproduction gives
\[
P_hE_hq=0.
\]
The injectivity of $E_h$ preserves the dimension of $\mathcal H_p$.
Since $E_h(\mathcal H_p)\subset\ker P_h$,
\[
\dim\ker P_h
\ge\dim E_h(\mathcal H_p)
=\dim\mathcal H_p
=2p+1.
\]

Taking into account, the full column rank of the stacked matrix means
\[
P_hc=0 \qandq B_hc=0
\quad\Longrightarrow\quad c=0.
\]
Thus, $B_h$ is injective in the space $\ker P_h$ of dimension $k_h$, 
which requires $m_\Gamma\ge k_h$.

Furthermore, if $P_h$ has full row rank, the rank--nullity gives
\[
k_h=n_h-n_-=n_{\rm rep}.
\]
Together with \eqref{eq:harmonic_nullity}, this proves
$n_{\rm rep}\ge2p+1$.
\end{proof}

For $p=4$, we obtain $\dim(\ker P_h)\ge 9$, which means that $\ker P_h$ may have a larger dimension. The polynomial reproduction follows from the exactness of the local interpolation and the identity $\sum_s w_s=1$ \cite{barrera2026tensor,dellaccio2016triangular}.

\subsection{Algebraic observability and square completion}

\begin{proposition}
\label{thm:observability}
\label{cor:square_err}
Let $P\in\mathbb{R}^{m_P\times n}$ and $B\in\mathbb{R}^{m_B\times n}$. Let $Z\in\mathbb{R}^{n\times k}$ have columns that form a basis of $\ker(P)$, where $k=n-\operatorname{rank}(P)$. Then
\begin{equation}
\rank\begin{bmatrix}P\\B\end{bmatrix}-\rank P=\rank(BZ).
\label{eq:observability}
\end{equation}
Consequently, the stacked matrix has full column rank if and only if
\begin{equation}
\rank(BZ)=k.
\label{eq:full_observability}
\end{equation}
If, in addition, $P$ has full row rank and $B$ has $k = n - m_P$ rows, the resulting stacked square matrix is nonsingular.
\end{proposition}

\begin{proof}
Define $A=[P^T\ B^T]^T$. Each $x\in\ker P$ admits a unique representation $x=Z\alpha$, and $x\in\ker A$ if and only if $BZ\alpha=0$. Therefore,
\[
\dim\ker A=\dim\ker(BZ)=k-\rank(BZ).
\]
Applying the rank–nullity theorem gives \eqref{eq:observability} and
\eqref{eq:full_observability}. Hence, $m_P+k=n$ implies that the matrix $A$ is square and has full column rank which means that $A$ is nonsingular.
\end{proof}

The theorem explains the different roles of the boundary completion methods. The methods ERR-P and ERR-U select boundary rows without using $Z_h$ and may leave unresolved modes. Also, the ERR-OS method adds boundary equations, but the full column rank of $A_h$ still requires \eqref{eq:full_observability}. However, in the ERR-QR method, we use this restriction to select the rows to complete the square system. We do not select boundary rows based on the uniform distribution of the added boundary point or the geometry of $\Gamma$. Related nonsingularity results derived from sampling-based analyses have also been established for other collocation methods, including certain randomly sampled Kansa-type systems \cite{cavoretto2025kansa}.

\begin{proposition}
\label{prop:square_completion}
\label{prop:qr_completion}
Let $P\in\R^{m\times n}$ have full row rank, set $k=n-m$, and let the columns of $Z\in\R^{n\times k}$ form a basis of $\ker P$. Let $C\in\R^{M\times n}$ with $M\ge k$. One can select from $C$ exactly $k$ rows that form a block $B$ such that
\[
A=\begin{bmatrix}P\\B\end{bmatrix}
\]
is nonsingular if and only if $\rank(CZ)=k$.

When $\rank(CZ)=k$, the column-pivoted QR applied to $(CZ)^T$ selects such a set in its first $k$ pivots. At each step, we assume that QR selects the remaining column with the largest residual norm.
\end{proposition}

\begin{proof}
If $B$ gives a nonsingular matrix $A$, then $\rank(BZ)=k$ by Proposition~\ref{thm:observability}. Since the rows of $BZ$ are rows of $CZ$, so $\rank(CZ)=k$. 

Conversely, if $\rank(CZ)=k$, we select $k$ linearly independent rows of $CZ$ to form the block $B$. Then $BZ$ is nonsingular, and Proposition~\ref{thm:observability} proves that the matrix $A$ is nonsingular. 

We set $D=(CZ)^T$, such that $\rank(D)=k$. We suppose also that the first $j-1$ selected columns are linearly independent, where $1\le j\le k$. If all remaining columns had zero residual after orthogonal projection onto the span of these selected columns, then every column of $D$ would belong to that span. This would imply
\[
\operatorname{rank}(D)\le j-1<k,
\]
which is a contradiction. Thus, at least one remaining residual has positive norm, and the pivoting rule selects a column with positive residual norm. By induction, the first $k$ pivot columns are linearly independent.

Equivalently, write
\[
D\Pi
=
Q\begin{bmatrix}R_{11}&R_{12}\end{bmatrix},
\qquad
R_{11}\in\mathbb{R}^{k\times k},
\]
the upper triangular matrix $R_{11}$ has nonzero diagonal entries and is therefore nonsingular. Since $J$ follows the pivot order,
\[
(BZ)^T=QR_{11}.
\]
Hence $BZ$ is nonsingular. Proposition~\ref{cor:square_err} then gives the nonsingularity of the stacked square matrix. If $k=0$, then $P$ is already a square matrix and nonsingular without adding boundary rows.

\end{proof}

The selection of the boundary rows from the candidate rows of $C$ uses the procedure of the classical column-pivoted Householder QR \cite{businger1965householder}. When fitted boundary rows are fixed in advance, we give the additional completion conditions. We first form
\[
P_F=
\begin{bmatrix}
P\\F
\end{bmatrix},
\]
where $P\in\R^{m\times n}$ and $F\in\R^{q\times n}$. We assume that $P_F$ has full row rank, and set
\[
k_F=n-m-q.
\]
Let the columns of $Z_F\in\R^{n\times k_F}$ form a basis of $\ker P_F$, and let $C\in\R^{M\times n}$ contain the unfitted boundary candidates, with $M\ge k_F$.

Then, we apply Proposition~\ref{prop:square_completion} to $P_F$ and find that exactly $k_F$ candidate rows can be selected to give a nonsingular square matrix
\[
\begin{bmatrix}
P\\F\\C_{J,:}
\end{bmatrix}
\]
if and only if
\[
\rank(CZ_F)=k_F.
\]
Here, the set $J$ contains the indices of the selected unfitted row. Under this condition, we use column-pivoted QR applied to $(CZ_F)^T$ to select the block $C_{J,:}$.

\subsection{Quantitative observability}
\begin{proposition}
\label{prop:quantitative_observability}
Let $\bar P\in\R^{m\times n}$ have rank $r$, with $0<r<n$.
Let $Z$ and $Y$ have orthonormal columns spanning
$\ker\bar P$ and its orthogonal complement, respectively.
Let $\bar B\in\R^{s\times n}$, with $s\ge n-r$, and set
\[
\alpha=\sigma_{\min}(\bar P Y),\qquad
\gamma=\sigma_{\min}(\bar B Z),\qquad
\beta=\|\bar B Y\|_2.
\]
Assume that $\gamma>0$. Let
\[
A=\begin{bmatrix}\bar P\\\bar B\end{bmatrix},
\qquad
K=\begin{bmatrix}
\alpha^{-1}&0\\
\beta/(\alpha\gamma)&\gamma^{-1}
\end{bmatrix}.
\]
Then $A$ has full column rank and
\begin{equation}
\|K\|_2^{-1}\le\sigma_{\min}(A)\le\gamma.
\label{eq:quantitative_observability}
\end{equation}
\end{proposition}
\begin{proof}
The restriction of $\bar P$ to $(\ker\bar P)^\perp$ is
injective, so $\alpha>0$.

We write $c=Ya+Zz$. Orthogonality gives
\[
\|c\|_2^2=\|a\|_2^2+\|z\|_2^2.
\]
Since $\bar PZ=0$,
\[
\|a\|_2\le\alpha^{-1}\|\bar Pc\|_2.
\]
Moreover,
\[
\gamma\|z\|_2
\le\|\bar BZz\|_2
\le\|\bar Bc\|_2+\beta\|a\|_2.
\]
Combining these inequalities gives the componentwise bound
\[
\begin{bmatrix}\|a\|_2\\\|z\|_2\end{bmatrix}
\le
K\begin{bmatrix}\|\bar Pc\|_2\\\|\bar Bc\|_2\end{bmatrix}.
\]
Taking Euclidean norms yields
\[
\|c\|_2\le\|K\|_2\|Ac\|_2.
\]
Thus $A$ has full column rank and
$\sigma_{\min}(A)\ge\|K\|_2^{-1}$.

For the upper bound, we choose a unit vector $z$ such that
$\|\bar BZz\|_2=\gamma$ and set $c=Zz$.
Then $\|c\|_2=1$ and $\bar Pc=0$, so
\[
\sigma_{\min}(A)\le\|Ac\|_2
=\|\bar BZz\|_2=\gamma.
\]
\end{proof}

For ERR-OS, we take
\[
\bar P=n_-^{-1/2}D_PP_h \qandq
\bar B=m_\Gamma^{-1/2}D_BB_h.
\]
The boundary margin $\gamma$ in Proposition~\ref{prop:quantitative_observability}  measures the boundary action on $\ker P_h$. The lower bound for $\sigma_{\min}(\widetilde A_h)$ also depends on the interior action $\alpha$ and the coupling term $\beta$.

If the two normalized and scaled blocks of $\widetilde A_h$ in \eqref{eq:weighted_ls} are not empty, then every row of $D_PP_h$ and $D_BB_h$ has unit Euclidean norm. Consequently,
\begin{equation}
\|\widetilde A_h\|_F^2
=\frac1{n_-}\|D_PP_h\|_F^2+
\frac1{m_\Gamma}\|D_BB_h\|_F^2=2.
\label{eq:frobenius_scaling}
\end{equation}
For full column rank, the sum of the $n_h$ squared singular values equals two. Hence
\[
\sigma_{\min}(\widetilde A_h)\le\sqrt{2/n_h}
\qandq
\|\widetilde A_h^\dagger\|_2\ge\sqrt{n_h/2}.
\]
This is a consequence of the normalization and does not indicate an unstable discretization.

\section{Numerical experiments}
\label{sec:experiments}

In this numerical section, we present six numerical experiments. Unless otherwise stated, all numerical experiments are conducted using the baseline configuration summarized in Table~\ref{tab:numerical_parameters}. Each experiment includes, as needed, adjustments to the parameters identified as the primary variables of interest in the study.

\begin{table}[h]
\centering
\caption{Common configuration used in the experiments.}
\label{tab:numerical_parameters}
\begin{tabular}{lll}
\toprule
Parameter & Symbol & Value \\
\midrule
Local subgrid size
    & $(p+1)\times(p+1)$ & $5\times5$ ($p=4$) \\
Subgrid step (\cite{elharrak2026gbmsc})
    & -- & $(1,1)$ \\
Shepard exponent
    & $\mu$ & $4$ \\
Interior sign tolerance
    & $\tau_{\rm in}$ & $10^{-14}$ \\
Retention threshold
    & $\vartheta_{\min}$ & $0.70$ \\
Boundary retraction tolerance
    & $|\phi(y)|$ & $\le 10^{-13}$ \\
Refinement levels in each direction
    & $N$ & $27,\,37,\,47,\,57$ \\
\bottomrule
\end{tabular}
\end{table}

The manufactured solution for the Poisson problem is as follows:
\begin{equation*}
u(x,y)=\sin(\pi x)\sin(\pi y),
\end{equation*}
with the associated source term given by
\begin{equation*}
f(x,y)=2\pi^2\sin(\pi x)\sin(\pi y).
\end{equation*}

For each experiment, we evaluate errors independently of the collocation points. The interior error is evaluated at the interior nodes of a $101 \times 101$ Cartesian grid within each domain, whereas the boundary error is evaluated using 1000 uniformly distributed points along the boundary. We report the root-mean-square error $e_{\RMS}$ and the maximum absolute error $ e_{\max}$, together with the boundary errors $e_{\Gamma,\max}$ and $e_{\Gamma,\RMS}$, and the RMS gradient error for the support and boundary-representation tests.

For spectral diagnostics, we use the scaled matrix
\begin{equation*}
\widetilde A_h
=
\begin{bmatrix}
n_-^{-1/2}D_PP_h\\
m_\Gamma^{-1/2}D_BB_h
\end{bmatrix}.
\end{equation*}
For a nonsingular square system, this scaling does not modify the exact solution, provided that the right-hand side is scaled in a consistent manner. In the ERR-OS formulation, it is incorporated into the least-squares objective \eqref{eq:weighted_ls}. Unless stated otherwise, the numerical rank is reported using a matrix-dependent tolerance.
\begin{equation}
\tau_{\rm rank}
=
\epsmach\max(m,n)\sigma_{\max}
\label{eq:rank_tolerance}
\end{equation}
for the matrix under consideration for ranking. $\epsmach$ is the machine precision of double-precision arithmetic. So $\operatorname{rank}_{\tau}(A)=\#\{i:\sigma_i(A)>\tau_{\rm rank}\}$. For the scaled matrix, we report the smallest singular value $\sigma_{\min}$ and spectral condition number $\kappaTwo$:
\begin{equation}
\sigmin(\widetilde A_h),
\qandq
\kappaTwo(\widetilde A_h)
=\sigma_{\max}(\widetilde A_h)/\sigma_{\min}(\widetilde A_h).
\label{eq:spectral_diagnostics}
\end{equation}

For every boundary completion method, we report at least the following algebraic quantities
\begin{equation}
\rank(P_h),\qquad k_h=n_h-\rank(P_h),\qquad
\rank(B_hZ_h),\qquad \rank(A_h).
\label{eq:rank_diagnostics}
\end{equation}
The tables report the diagnostics available for each experiment. Computed nullities and rank classifications use numerical thresholds and must be distinguished from the exact quantities in Section~\ref{sec:analysis}. A matrix-dependent threshold does not need to give rank counts that satisfy the exact rank identity across independently scaled matrices. Under these tests, the boundary restriction can be classified as full rank while the complete stack is classified as deficient. This is how the annular ERR-U entry in Table~\ref{tab:baseline_refinement} should be interpreted.

For an orthonormal nullspace basis $Z_h$, we additionally report the boundary margin
\[
\gamma_h=\sigma_{\min}(D_BB_hZ_h).
\]
$\gamma_h$ measures the boundary action on $\ker P_h$ after normalization of the rows. The averaged margin $\tilde\gamma_h$ for the boundary block in $\widetilde A_h$ is $\gamma_h/\sqrt{m_\Gamma}$. If we replace each boundary row with $q$ identical rows, we obtain the new value $\gamma_h$ by rescaling $\sqrt q$, while the averaged margin $\tilde\gamma_h$ remains unchanged. Therefore, comparisons of $\gamma_h$ values in experiments with different numbers of rows should explicitly account for the dependence of $\gamma_h$ on the sample size. We retain the original reported values. The margin alone does not control the stability of the full system; see Proposition~\ref{prop:quantitative_observability}.

In each run, we use the labels \emph{full column rank}, \emph{deficient}, or \emph{invalid budget}. Deficient cases are not included in convergence-rate fits. An invalid budget case is used when the specified number of rows does not provide the required boundary equations for the intended square completion.

\subsection{Smooth implicit-domain refinement}
\label{sec:baseline_refinement}

This experiment characterizes the baseline performance of the four boundary completion methods on smooth implicit domains. In this study, we consider a disk and an annulus represented by exact signed-distance fields, as well as a smooth star domain described by a level-set function. All computational domains are embedded within the unit square. Let
\[
c=(0.5,0.5),\qquad
r(x)=\|x-c\|_2,
\qquad
\theta(x)=\operatorname{atan2}(x_2-0.5,x_1-0.5).
\]
The three smooth implicit geometries used in the main numerical experiments are    
\[
\begin{aligned}
\Omega_{\rm disk}
&=\{x:r(x)<0.40\},\\
\Omega_{\rm ann}
&=\{x:0.15<r(x)<0.43\},\\
\Omega_{\rm star}
&=\{c\}\cup\left\{x\ne c:
r(x)<0.34+0.08\cos\!\bigl(5\theta(x)\bigr)
\right\}.
\end{aligned}
\]

Accordingly, the disk and the annulus are represented by the exact signed-distance
functions
\[
\phi_{\rm disk}(x)=r(x)-0.40 \qandq
\phi_{\rm ann}(x)
=
\max\{r(x)-0.43,\;0.15-r(x)\},
\]
and the star is represented by the regular level-set function
\[
\phi_{\rm star}(x)
=
r(x)-0.34-0.08\cos\!\bigl(5\theta(x)\bigr),
\]
which is not an exact signed-distance function. $\phi_{\rm star}$ is also a regular function near the boundary. In the center $c$, the polar representation is not defined; therefore, we consider explicitly the center $c$ as an interior point. No global continuity of this polar sign field is required.
   
ERR-P and ERR-U are the square geometric baselines. ERR-OS uses boundary oversampling with $\rho=1.5$, and ERR-QR uses the same candidate-cloud construction. The finest-grid table reports the replacement budget $n_{\rm rep}$, the nullity of numerical PDE $k_h$, the boundary rank $\rank(B_hZ_h)$, the status of the complete system, the scaled spectral diagnostics, and the RMS error $e_{\rm RMS}$.
   
\begin{figure} [h]
\centering
\begin{subfigure}[t]{0.32\textwidth}
\centering
\includegraphics[width=\textwidth]{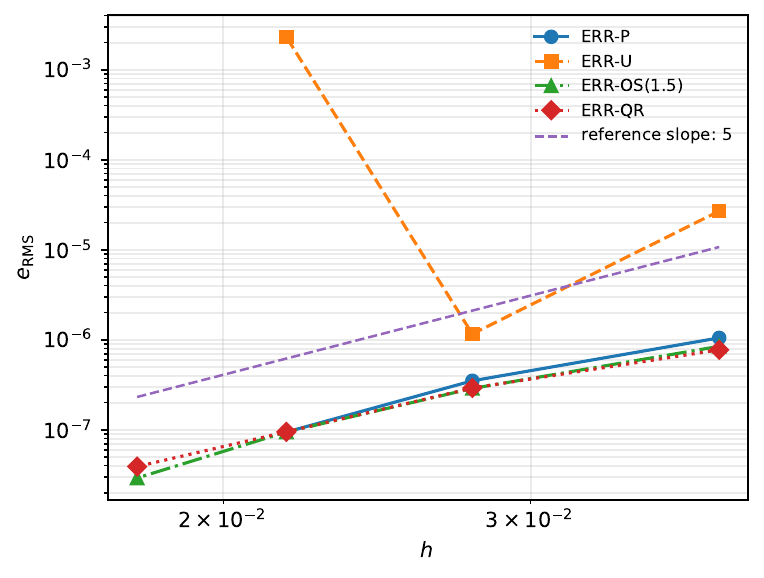}
\caption{Disk.}
\end{subfigure}\hfill
\begin{subfigure}[t]{0.32\textwidth}
\centering
\includegraphics[width=\textwidth]{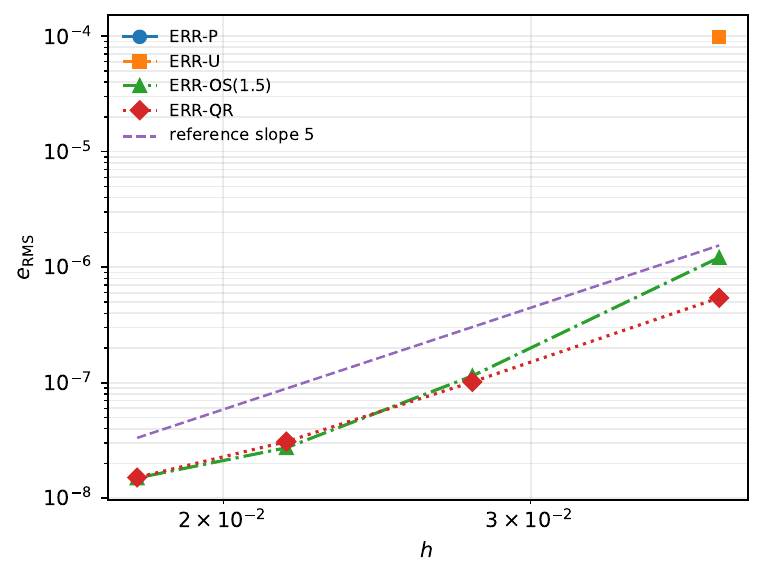}
\caption{Annulus.}
\end{subfigure}\hfill
\begin{subfigure}[t]{0.32\textwidth}
\centering
\includegraphics[width=\textwidth]{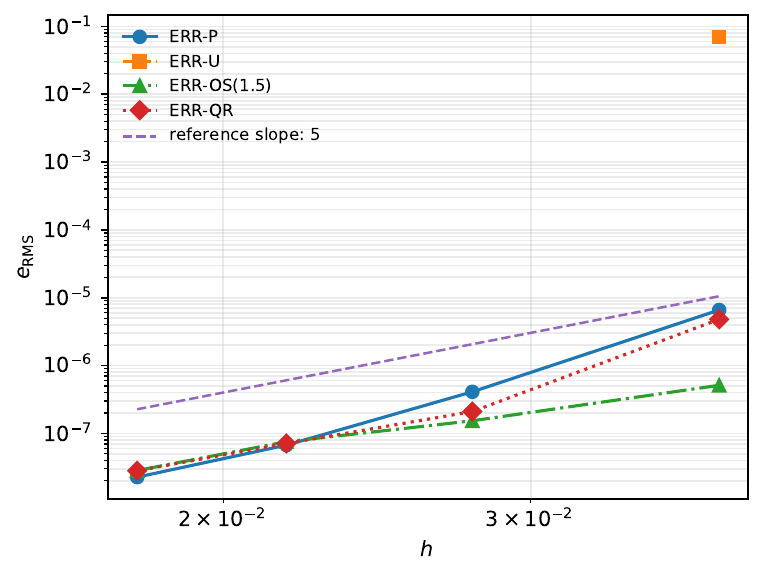}
\caption{Star.}
\end{subfigure}
\caption{Log-log plots of the error $e_{\RMS}$ as a function of the grid spacing $h$ for the four completion methods (ERR-P, ERR-U, ERR-OS(1.5), and ERR-QR): (a) disk, (b) annulus, and (c) star. RMS value errors are plotted only for solved cases.}
\label{fig:baseline_refinement}
\end{figure}

Table \ref{tab:baseline_refinement} presents a comparison with $N=57$ for smooth implicit-domain geometries (disk, annulus, and star) using four completion methods (ERR-P, ERR-U, ERR-OS(1.5) and ERR-QR). Figure \ref{fig:baseline_refinement} shows log--log plots of $e_{\RMS}$, as a function of grid spacing $h$ for the four completion methods in the following geometries: disk (Figure~\ref{fig:baseline_refinement}-(a)), annulus (Figure~\ref{fig:baseline_refinement}-(b)), and star (Figure~\ref{fig:baseline_refinement}-(c)). The RMS errors are plotted only for full-rank cases. For the ERR-OS approach, the parameter $\rho$ is fixed at 1.5. The square completion methods ERR-P and ERR-U are not uniformly admissible in this refinement set. ERR-P is full column rank in 8 of the 12 domain--grid cases and ERR-U in 5 of 12. ERR-OS$(1.5)$ and ERR-QR are full column rank in all 12 cases. At $N=57$, their RMS errors range from $1.49e{-8}$ to $3.92e{-8}$. The convergence slopes reported in Figure~\ref{fig:baseline_refinement} indicate that ERR-OS(1.5) and ERR-QR achieve the fifth-order accuracy across the considered levels of grid refinement.

\begin{table} [h]
\centering
\caption{Report of finest-grid comparison at $N=57$ for smooth implicit-domain geometries (the disk, annulus, and star) using the four completion methods (ERR-P, ERR-U, ERR-OS(1.5), and ERR-QR). For ERR-OS, the parameter $\rho$ is set to 1.5.}
\label{tab:baseline_refinement}
\small
 
\begin{tabular}{llllllllll}
\toprule
Domain & completion method & $n_h$ & $n_{\rm rep}$ & $k_h$ & $\operatorname{rank}(B_hZ_h)$ & Status & $\sigma_{\min}(\tilde A_h)$ & $\kappa_2(\tilde A_h)$ & $e_{\rm RMS}$  \\
\midrule
\multirow{4}{*}{Disk} & ERR-P & 1809 & 228 & 228 & 226 & deficient & 2.57e-18 & -- & --  \\
 & ERR-U & 1809 & 228 & 228 & 227 & deficient & 3.67e-18 & -- & --  \\
 & ERR-OS$(1.5)$ & 1809 & 228 & 228 & 228 & full & 4.68e-06 & 2.28e+04 & 2.92e-08  \\
 & ERR-QR & 1809 & 228 & 228 & 228 & full & 1.04e-05 & 1.12e+04 & 3.92e-08  \\
\midrule
\multirow{4}{*}{Annulus} & ERR-P & 1936 & 344 & 344 & 344 & full & 4.63e-06 & 1.95e+04 & 1.49e-08  \\
 & ERR-U & 1936 & 344 & 344 & 344 & deficient & 1.08e-14 & -- & --  \\
 & ERR-OS$(1.5)$ & 1936 & 344 & 344 & 344 & full & 2.82e-05 & 3.33e+03 & 1.49e-08  \\
 & ERR-QR & 1936 & 344 & 344 & 344 & full & 4.15e-05 & 2.22e+03 & 1.51e-08  \\
\midrule
\multirow{4}{*}{Star} & ERR-P & 1409 & 237 & 237 & 237 & full & 3.46e-07 & 3.30e+05 & 2.26e-08  \\
 & ERR-U & 1409 & 237 & 237 & 233 & deficient & 2.91e-19 & -- & --  \\
 & ERR-OS$(1.5)$ & 1409 & 237 & 237 & 237 & full & 1.56e-05 & 6.96e+03 & 2.84e-08  \\
 & ERR-QR & 1409 & 237 & 237 & 237 & full & 1.97e-05 & 5.67e+03 & 2.82e-08  \\
\bottomrule
\end{tabular}
 
\end{table}

\subsection{Effect of the oversampling ratio}
\label{sec:rho_exp}

We study ERR-OS at $N=47$ using
\begin{equation}
\rho\in\{1.10,1.25,1.50,1.75,2.00\}.
\label{eq:rho_set}
\end{equation}
The objective of \eqref{eq:rho_set} is not to optimize $\rho$, but rather to assess whether the behavior observed with $\rho = 1.5$ remains robust under moderate oversampling. Consequently, we vary $\rho$ in accordance with \eqref{eq:rho_set} while keeping fixed the computational domain, the retained source space, the candidate cloud, the row normalization procedure, and the least-squares weighting scheme. This experiment tests whether a modest amount of boundary oversampling is enough to remove rank failures and whether increasing the oversampling materially changes the error or the spectral diagnostics. 

\begin{figure} [h]
\centering
\begin{subfigure}[t]{0.4\textwidth}
\centering
\includegraphics[width=\textwidth]{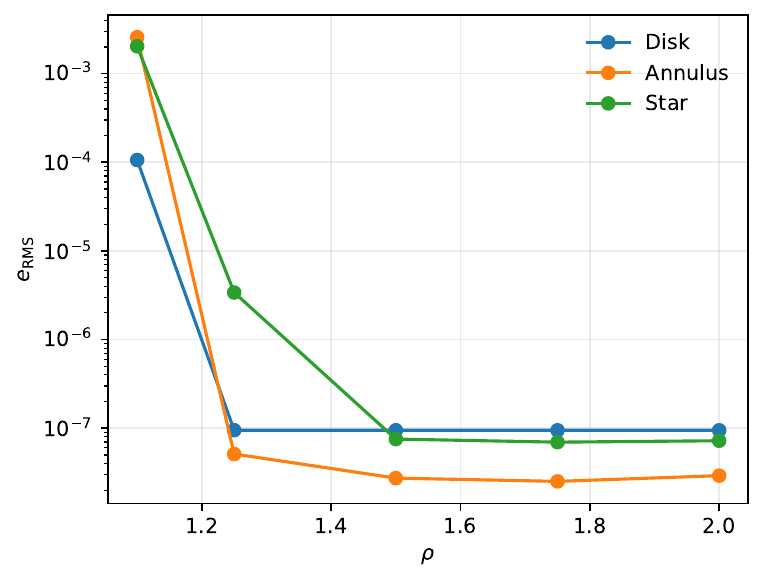}
\caption{$e_{\rm RMS}$.}
\end{subfigure}
\begin{subfigure}[t]{0.4\textwidth}
\centering
\includegraphics[width=\textwidth]{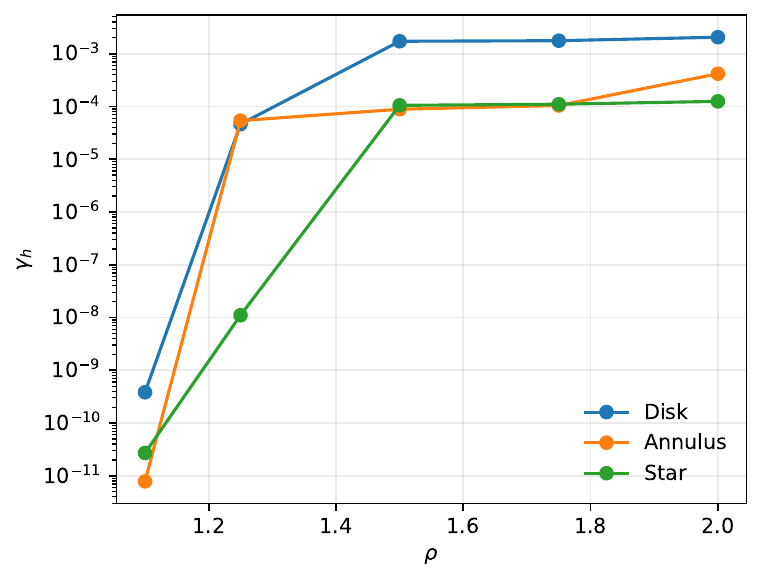}
\caption{$\gamma_h$.}
\end{subfigure}
\caption{Sensitivity analysis of ERR-OS with respect to the boundary oversampling ratio $\rho$ for $N=47$. The approximation error $e_{\rm RMS}$ and algebraic observability diagnostics $\gamma_h$ are reported in panels (a) and (b), respectively.}
\label{fig:rho_study}
\end{figure}

Table \ref{tab:rho_study} reports an ERR-OS oversampling study conducted at $N=47$ for several values of the oversampling ratio, $\rho \in \{1.10, 1.25, 1.50, 1.75, 2.00\}$, with performance assessed separately within each domain. Figure \ref{fig:rho_study} presents a sensitivity analysis of ERR-OS with respect to the ratio $\rho$ of the boundary oversampling for $N=47$. The approximation error $e_{\rm RMS}$ and algebraic observability diagnostics $\gamma_h$ are reported in Figure~\ref{fig:rho_study}-(a) and Figure~\ref{fig:rho_study}-(b), respectively.

\begin{table} [h]
\centering
\caption{ERR-OS oversampling study at $N=47$ for multiple values of $\rho$, $\rho\in\{1.10,1.25,1.50,1.75,2.00\}$, evaluated within each domain. Within each domain, the source space, PDE block, ordered candidate cloud, and normalization rule are fixed.}
\label{tab:rho_study}
\small
 
\begin{tabular}{llllllll}
\toprule
Domain & $\rho$ & $m_\Gamma$ & $\operatorname{rank}(B_hZ_h)$ & Status & $\gamma_h$ & $\sigma_{\min}(\widetilde A_h)$ & $e_{\rm RMS}$ \\
\midrule
\multirow{5}{*}{Disk}   & 1.10 & 212 & 192 & full & 3.83e-10 & 2.62e-11 & 1.06e-04 \\
 & 1.25 & 240 & 192 & full & 4.63e-05 & 2.86e-06 & 9.49e-08 \\
 & 1.50 & 288 & 192 & full & 1.72e-03 & 8.91e-05 & 9.48e-08 \\
 & 1.75 & 336 & 192 & full & 1.77e-03 & 8.90e-05 & 9.47e-08 \\
 & 2.00 & 384 & 192 & full & 2.06e-03 & 8.91e-05 & 9.47e-08 \\
\midrule
\multirow{5}{*}{Annulus}   & 1.10 & 304 & 276 & full & 7.84e-12 & 4.48e-13 & 2.57e-03 \\
 & 1.25 & 345 & 276 & full & 5.39e-05 & 2.74e-06 & 5.11e-08 \\
 & 1.50 & 414 & 276 & full & 8.87e-05 & 4.28e-06 & 2.73e-08 \\
 & 1.75 & 483 & 276 & full & 1.04e-04 & 4.69e-06 & 2.51e-08 \\
 & 2.00 & 552 & 276 & full & 4.18e-04 & 1.72e-05 & 2.91e-08 \\
\midrule
\multirow{5}{*}{Star}  & 1.10 & 226 & 205 & full & 2.71e-11 & 1.77e-12 & 2.03e-03 \\
  & 1.25 & 257 & 205 & full & 1.10e-08 & 6.75e-10 & 3.39e-06 \\
  & 1.50 & 308 & 205 & full & 1.05e-04 & 5.55e-06 & 7.54e-08 \\
  & 1.75 & 359 & 205 & full & 1.10e-04 & 5.68e-06 & 6.96e-08 \\
  & 2.00 & 410 & 205 & full & 1.25e-04 & 5.97e-06 & 7.22e-08 \\
\bottomrule
\end{tabular}
 
\end{table}

The numerical rank tests classify all fifteen systems in Table~\ref{tab:rho_study} as full column rank. At $\rho=1.10$, the smallest singular values are very small, and the RMS errors range from $1.06e{-4}$ to $2.57e{-3}$. Increasing $\rho$ to $1.25$ reduces the disk and annulus errors to $9.49e{-8}$ and $5.11e{-8}$, respectively, while the star error is $3.39e{-6}$. For $1.50\le\rho\le2.00$, the errors on all three domains are of order $10^{-8}$ to $10^{-7}$ at the reported precision. Numerical full rank alone therefore does not explain the differences in accuracy. These results support $\rho=1.5$ as a reference choice for the tested configurations, but they establish neither a universal threshold nor an optimal ratio. Comparisons of the raw boundary margin must also account for its dependence on the number of boundary rows.

\subsection{Comparison between ERR-OS and ERR-QR}
\label{sec:os_qr_exp}

We use an identical candidate cloud and the same candidate matrix $C_h$ for both completion methods. ERR-OS chooses a quasi-uniform subset of cardinality $\lceil\rho n_{\rm rep}\rceil$ rows. Under the boundary rank condition, ERR-QR selects $k_h$ rows such that their restrictions to $\ker P_h$ are linearly independent. The ERR-QR selection need not be contained in the ERR-OS subset. The comparison therefore concerns two selections from one common candidate family, with additional observations retained in ERR-OS and an algebraic completion selected in ERR-QR.
   
For each domain, we record the candidate count $M_h$, the boundary-row count $m_\Gamma$, $\rank(B_hZ_h)$, $\gamma_h$, $\sigmin(\widetilde A_h)$, the RMS error $e_{\mathrm{RMS}}$, and the boundary error $e_{\Gamma,\max}$. Table~\ref{tab:os_qr} presents
a comparison between ERR-OS(1.5) and ERR-QR using the same oversampled candidate points for $N=47$. Both boundary completions in each domain use the same candidate matrix.

\begin{table} [h]
\centering
\caption{Comparison between ERR-OS(1.5) and ERR-QR using the same oversampled candidate points for $N=47$. Both boundary completion methods in each domain use the same candidate matrix.}
\label{tab:os_qr}
\small
 
\begin{tabular}{llllllllll}
\toprule
Domain & completion method & $M_h$ & $m_\Gamma$ & $\mathrm{rank}(B_hZ_h)$ & $\gamma_h$ & Status & $\sigma_{\min}(\widetilde A_h)$ & $e_{\mathrm{RMS}}$ & $e_{\Gamma,\max}$ \\
\midrule
\multirow{2}{*}{Disk}  & ERR-OS(1.5) & 1803 & 288 & 192 & 1.72e-03 & full & 8.91e-05 & 9.48e-08 & 9.37e-10 \\
 & ERR-QR & 1803 & 192 & 192 & 1.75e-03 & full & 8.91e-05 & 9.51e-08 & 1.25e-09 \\
 \midrule
\multirow{2}{*}{Annulus} & ERR-OS(1.5) & 2606 & 414 & 276 & 8.87e-05 & full & 4.28e-06 & 2.73e-08 & 1.05e-08 \\
 & ERR-QR & 2606 & 276 & 276 & 3.11e-04 & full & 1.76e-05 & 3.10e-08 & 3.26e-09 \\
 \midrule
\multirow{2}{*}{Star} & ERR-OS(1.5) & 2090 & 308 & 205 & 1.05e-04 & full & 5.55e-06 & 7.54e-08 & 3.53e-08 \\
 & ERR-QR & 2090 & 205 & 205 & 2.50e-04 & full & 1.32e-05 & 7.07e-08 & 3.38e-08 \\
\bottomrule
\end{tabular}
 
\end{table}

Both completion methods are classified as full column rank in the three paired cases. ERR-QR uses $n_{\rm rep}$ boundary rows, whereas ERR-OS$(1.5)$ uses $\lceil1.5n_{\rm rep}\rceil$; QR therefore removes approximately one third of the boundary equations. The RMS errors are comparable: $9.48e{-8}$ compared to $9.51e{-8}$ on the disk, $2.73e{-8}$ compared to $3.10e{-8}$ on the annulus, and $7.54e{-8}$ compared to $7.07e{-8}$ on the star, for ERR-OS and ERR-QR, respectively. On the annulus and star, ERR-QR has larger reported boundary margins and larger scaled smallest singular values. The latter are not condition numbers. These cases show that compression can retain comparable sampled accuracy with fewer boundary equations. They do not establish a universal preference between the completion methods or a reduction in total computational cost, since the nullspace and selection costs are not reported.

\subsection{Trimming and algebraic admissibility}
\label{sec:trimming_exp}

Trimming changes both the retained source space and the available replacement budget. We therefore vary $\vartheta_{\min}$ and the reference sequence is
\begin{equation}
\vartheta_{\min}\in\{0.3, 0.4, 0.5, 0.6, 0.7, 0.8, 0.9 ,0.96\}.
\label{eq:trim_values}
\end{equation}
For each value of $\vartheta_{\min}$, both completion methods employ an identical source space and the same interior PDE block. In this experiment, we record $n_h$, $n_{\rm rep}$, $\dim\ker P_h$, $\rank(B_hZ_h)$, and the status of the complete-system using $N=37$. This constitutes a direct numerical evaluation of the distinction between the replacement budget $n_{\rm rep}$ and PDE nullity $k_h$.  

\begin{figure} [h]
\centering
\begin{subfigure}[t]{0.47\textwidth}
\centering
\includegraphics[width=\textwidth]{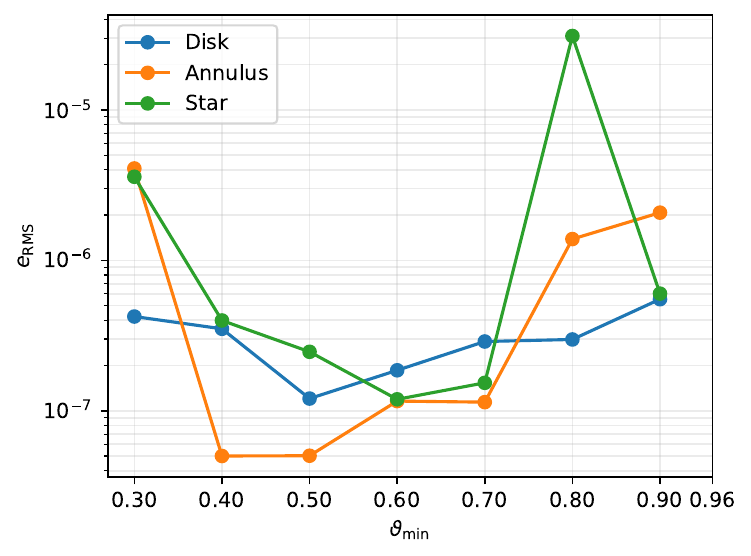}
\caption{RMS error for ERR-OS(1.5).}
\end{subfigure}
\begin{subfigure}[t]{0.47\textwidth}
\centering
\includegraphics[width=\textwidth]{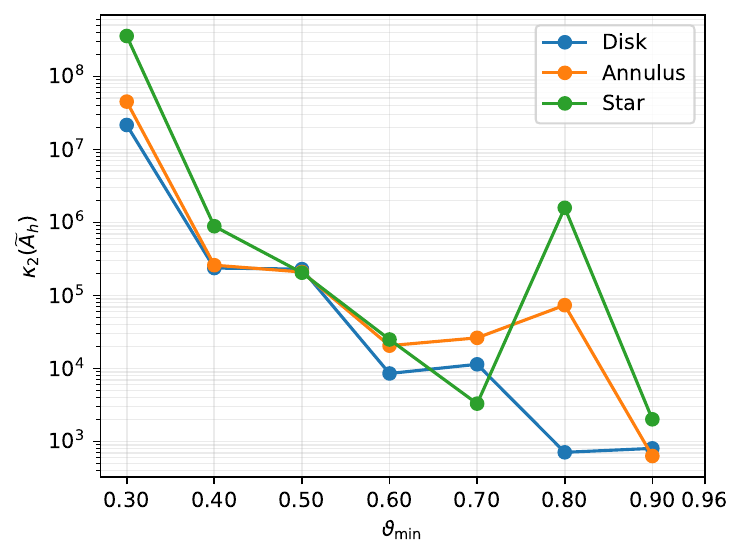}
\caption{$\kappa_2(\tilde A_h)$ for ERR-OS(1.5).}
\end{subfigure}\\
\begin{subfigure}[t]{0.47\textwidth}
\centering
\includegraphics[width=\textwidth]{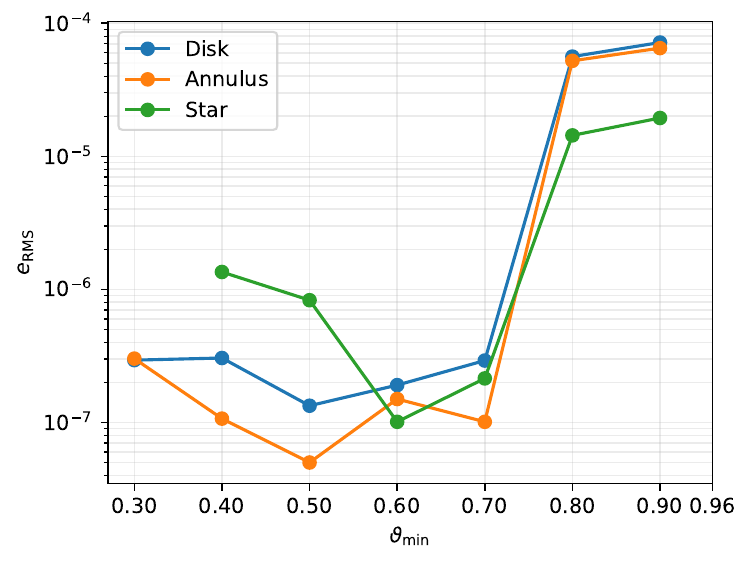}
\caption{RMS error for ERR-QR.}
\end{subfigure}
\begin{subfigure}[t]{0.47\textwidth}
\centering
\includegraphics[width=\textwidth]{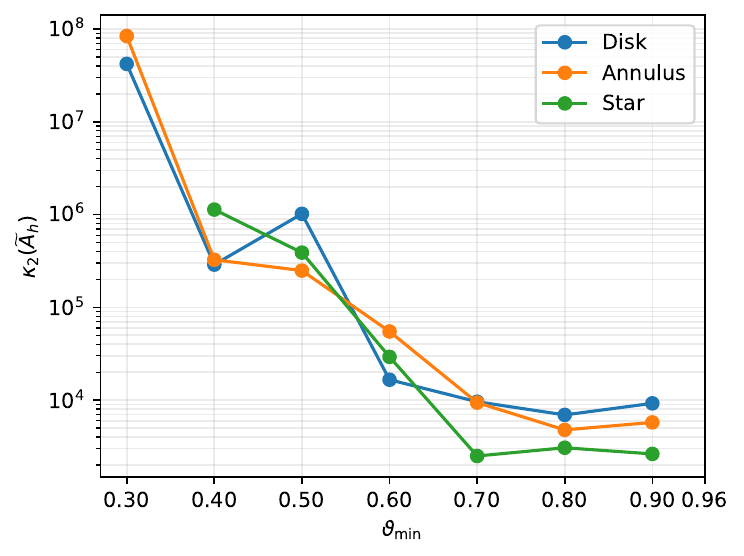}
\caption{$\kappa_2(\tilde A_h)$ for ERR-QR.}
\end{subfigure}
\caption{Semi-log plots of the RMS error and the spectral condition number $\kappaTwo$ as a function of the trimming threshold $\vartheta_{\min}$ using ERR-OS(1.5) and ERR-QR for the three smooth geometries at $N=37$. Left panels (a) and (c) present the RMS errors as a function of $\vartheta_{\min}$. Right panels (b) and (d) present the corresponding $\kappaTwo$.}
\label{fig:trimming_admissibility}
\end{figure}

Table~\ref{tab:trimming} summarizes the influence of $\vartheta_{\min}$ on ERR-OS(1.5) and ERR-QR for the disk. Figure~\ref{fig:trimming_admissibility} reports semi-logarithmic plots of $e_{\RMS}$ and $\kappaTwo$ as functions of $\vartheta_{\min}$, obtained using ERR-OS(1.5) and ERR-QR. The left panels (Figure~\ref{fig:trimming_admissibility}-(a) and Figure~\ref{fig:trimming_admissibility}-(c)) report the RMS error as a function of $\vartheta_{\min}$. The right panels (Figure~\ref{fig:trimming_admissibility}-(b) and Figure~\ref{fig:trimming_admissibility}-(d)) report the corresponding spectral condition number $\kappaTwo$.

\begin{table}[h]
\centering
\caption{Effect of the trimming threshold $\vartheta_{\min}$ on ERR-OS(1.5)
and ERR-QR for the disk at $N=37$.}
\label{tab:trimming}
\small
 
\begin{adjustbox}{max width=\textwidth}
\begin{tabular}{llllllllllll}
\toprule
& & & &
\multicolumn{4}{c}{ERR-OS} &
\multicolumn{4}{c}{ERR-QR} \\
\cmidrule(lr){5-8}
\cmidrule(lr){9-12}
$\vartheta_{\min}$
& $n_h$
& $n_{\rm rep}$
& $k_h$
& $\operatorname{rank}(B_hZ_h)$
& Status
& $\sigma_{\min}(\widetilde A_h)$
& $e_{\rm RMS}$
& $\operatorname{rank}(B_hZ_h)$
& Status
& $\sigma_{\min}(\widetilde A_h)$
& $e_{\rm RMS}$ \\
\midrule
0.30 & 997 & 340 & 340
& 340 & full & 6.16e-09 & 4.23e-07
& 340 & full & 3.13e-09 & 2.94e-07 \\

0.40 & 921 & 264 & 264
& 264 & full & 5.62e-07 & 3.51e-07
& 264 & full & 4.85e-07 & 3.05e-07 \\

0.50 & 897 & 240 & 240
& 240 & full & 5.90e-07 & 1.21e-07
& 240 & full & 1.44e-07 & 1.34e-07 \\

0.60 & 825 & 168 & 168
& 168 & full & 1.59e-05 & 1.86e-07
& 168 & full & 9.32e-06 & 1.91e-07 \\

0.70 & 801 & 144 & 144
& 144 & full & 1.23e-05 & 2.89e-07
& 144 & full & 1.57e-05 & 2.92e-07 \\

0.80 & 733 & 76 & 76
& 76 & full & 1.86e-04 & 2.98e-07
& 76 & full & 2.38e-05 & 5.62e-05 \\

0.90 & 705 & 48 & 48
& 48 & full & 1.84e-04 & 5.52e-07
& 48 & full & 1.99e-05 & 7.17e-05 \\

0.96 & 657 & 0 & 9
& -- & invalid budget & -- & --
& -- & invalid budget & -- & -- \\
\bottomrule
\end{tabular}
\end{adjustbox}
 
\end{table}

Table~\ref{tab:trimming} shows that increasing $\vartheta_{\min}$ reduces the active-source $n_h$ count and the replacement budget $n_{\rm rep}$. For $0.3\le\vartheta_{\min}\le0.9$, the numerical nullity equals $n_{\rm rep}$, and both completion methods are classified as full column rank. At $\vartheta_{\min}=0.96$, the strict retention rule with $p=4$ keeps only subgrids whose $25$ nodes satisfy $\phi<0$. Then $\I_+$ is empty, although nodes in the class $\I_0$ may remain. In the case of disk, we have an empty replacement budget and the numerical nullity $k_h$ is nine. The value $k_h=9$ is consistent with the bound $k_h\ge9$ in Proposition~\ref{prop:polynomial_kernel}, which does not give an upper bound on the kernel dimension. Without boundary rows, the harmonic modes remain unconstrained.

Figure~\ref{fig:trimming_admissibility} shows that with moderate trimming, the two completion methods perform similarly. For \(0.3\leq\vartheta_{\min}\leq0.7\), both achieve RMS errors of order \(10^{-7}\) despite large changes in \(n_h\), \(n_{\rm rep}\), and the smallest singular value. With higher trimming (\(\vartheta_{\min}=0.8\) and \(0.9\)), ERR-OS keeps RMS errors below \(6e{-7}\), while ERR-QR increases to \(5.62e{-5}\) and \(7.17e{-5}\). This shows that full column rank alone does not ensure comparable accuracy under high trimming, and ERR-OS benefits from preserving extra boundary information by oversampling. This result sets a practical upper limit on trimming for this construction and shows that the geometric trimming criterion must align with the algebraic nullity of the interior PDE block.

\subsection{Support sensitivity and boundary locality}
\label{sec:support}

This experiment compares GF and NA under ERR-QR on the disk, annulus, and star. The source nodes, the interior PDE block, the candidate cloud, and the selected boundary collocation points are fixed. The boundary evaluations and reconstructed values use the two support policies. 

In fact, we use the NA policy to keep the off-source boundary rows sparse. The PDE rows are already sparse under GF at an active source node $x_i$, because the weights of subgrids that do not contain the node $x_i$ and their derivatives up to order two vanish when $\mu>2$. Consequently, the PDE rows that use the GF or NA policy coincide.  Dirichlet evaluation at an active source is also identical under both policies. The difference between the two policies concerns off-source boundary points, where GF may involve all active coefficients

This experiment examines the accuracy and conditioning consequences for the specified common configuration. In this experiment, using $N=57$, we report the interior errors $e_{\rm RMS}$ and $e_{\max}$, the boundary errors $e_{\Gamma,\rm RMS}$ and $e_{\Gamma,\max}$, and the spectral condition number $\kappa_2(\widetilde A_h)$. Table~\ref{tab:support} reports the sensitivity of ERR-QR to the GF and NA support policies.

\begin{table}[h]
\centering
\caption{ERR-QR Support sensitivity at $N=57$.}
\label{tab:support}
\small
\begin{tabular}{lllllll}
\toprule
Domain & Support
& $e_{\rm RMS}$
& $e_{\max}$
& $e_{\Gamma,\rm RMS}$
& $e_{\Gamma,\max}$
& $\kappa_2(\widetilde A_h)$ \\
\midrule
\multirow{2}{*}{Disk}
& NA & 3.92e-08 & 6.23e-08 & 2.83e-09 & 1.58e-08 & 1.12e+04 \\
& GF & 3.92e-08 & 6.23e-08 & 2.83e-09 & 1.58e-08 & 1.12e+04 \\
\midrule
\multirow{2}{*}{Annulus}
& NA & 1.51e-08 & 4.42e-08 & 1.78e-10 & 8.90e-10 & 2.23e+03 \\
& GF & 1.51e-08 & 4.42e-08 & 1.78e-10 & 8.90e-10 & 2.23e+03 \\
\midrule
\multirow{2}{*}{Star}
& NA & 2.82e-08 & 5.10e-08 & 3.21e-09 & 2.68e-08 & 5.67e+03 \\
& GF & 2.82e-08 & 5.10e-08 & 3.21e-09 & 2.68e-08 & 5.67e+03 \\
\bottomrule
\end{tabular}
\end{table} 

At $N=57$, Table~\ref{tab:support} shows identical GF and NA errors and condition numbers. The maximum entrywise differences between their boundary matrices are $1.97e{-12}$, $4.31e{-12}$, and $2.81e{-11}$ for the disk, annulus, and star, respectively. Thus, for these smooth test domains, the additional subgrids included by GF have a numerically negligible effect at the selected boundary collocation points. NA therefore provides the same observed accuracy and conditioning while retaining a local boundary sparse evaluation.

\FloatBarrier
\subsection{Implicit, hybrid, and fitted boundary representations}
\label{sec:boundary_hierarchy}
 
In the case of a fitted component, an active source located precisely on the boundary contributes one Dirichlet row and does not need retraction. An implicit component supplies off-source boundary rows through the selected completion method. A hybrid formulation combines the two types. The complete system is square for a square completion under its rank and row-count assumptions, and overdetermined for ERR-OS. In this section, we use the ERR-QR completion for the boundary treatment of implicit boundary components.

The square hole is
\[
\Omega_{\rm hole}=(0.10,0.90)^2\setminus
\overline{B_{0.13}(0.5,0.5)},
\]
where $B_{0.13}(0.5,0.5)$ is the open disk of radius $0.13$ centered at $(0.5,0.5)$. Its outer rectangular boundary is fitted and its circular inner boundary is implicit. The fitted domain is the L-shaped domain
\[
\Omega_L=
(0,1)^2\setminus([0.5,1]\times[0.5,1]),
\]
whose boundary is fully fitted to the grid.

\begin{figure} [h]
\centering
\includegraphics[width=0.98\textwidth]{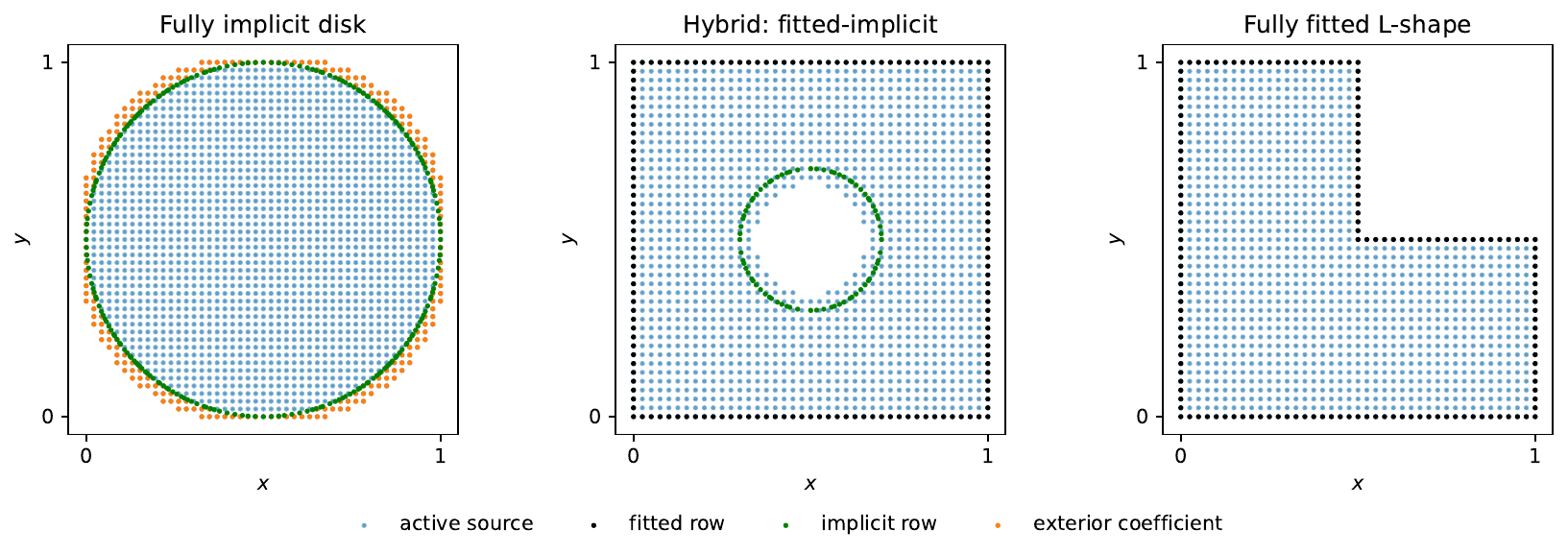}
\caption{Boundary structures for a fully implicit domain (disk), a hybrid fitted-implicit domain (square hole) and a fully fitted domain (L-shape).}
\label{fig:boundary_hierarchy}
\end{figure}

The algebraic reduction ensured by fitted components makes the row-selection algorithm for implicit boundary independent of the fitted boundary row. Figure~\ref{fig:boundary_hierarchy} shows the fully implicit, hybrid fitted--implicit, and fitted configurations. In the hybrid configuration, boundary conditions are applied independently on the outer rectangular boundary and on the inner circular boundary. In the third case, we use a fitted L-shaped domain to validate the fitted reduction on a non-rectangular geometry, as its boundary coincides with the source nodes.

Let $m_{\rm fit}$ denote the number of fitted Cartesian boundary rows and $m_{\rm imp}$ the number of rows obtained via the implicit-boundary completion. The total number of boundary rows is given by
\[
m_\Gamma=m_{\rm fit}+m_{\rm imp}.
\]
Table~\ref{tab:boundary_hierarchy} reports the fitted and implicit boundary-row counts, together with the corresponding geometric configuration.

\begin{table} [h]
\centering
\caption{Boundary-representation structure.}
\label{tab:boundary_hierarchy}
\small
 
\begin{tabular}{lrrrrrrrrlll}
\toprule
Case & $N$ & $n_h$ & $n_{\rm -}$ & $n_+$ & $m_{\rm fit}$ & $m_{\rm imp}$ & $m_\Gamma$ & Status & $e_{\mathrm{RMS}}$ & $e_{\Gamma,\max}$ & $e_{\Gamma,\mathrm{RMS}}$ \\
\midrule
Implicit disk & 47 & 1881 & 1649 & 228 & 0 & 232 & 232 & full & 8.28e-08 & 5.87e-09 & 1.06e-09 \\
Hybrid square--hole & 41 & 1560 & 1324 & 72 & 160 & 76 & 236 & full & 1.48e-07 & 5.44e-09 & 1.10e-09 \\
Fitted L--shape & 41 & 1281 & 1121 & 0 & 160 & 0 & 160 & full & 1.11e-07 & 8.94e-08 & 8.76e-09 \\
\bottomrule
\end{tabular}
\end{table} 

For the hybrid case, Table~\ref{tab:boundary_hierarchy} reports $m_\Gamma=236$, with $160$ outer-boundary rows and $76$ implicit inner-boundary rows. The RMS errors are $8.28e{-8}$ for the implicit disk, $1.48e{-7}$ for the hybrid case, and $1.11e{-7}$ for the fitted L-shaped domain. In the hybrid case, we first impose the fitted rows and then complete the square completion by implicit boundary rows. 

\begin{figure} [h]
\centering
\begin{subfigure}[t]{0.4\textwidth}
\centering
\includegraphics[width=\textwidth]{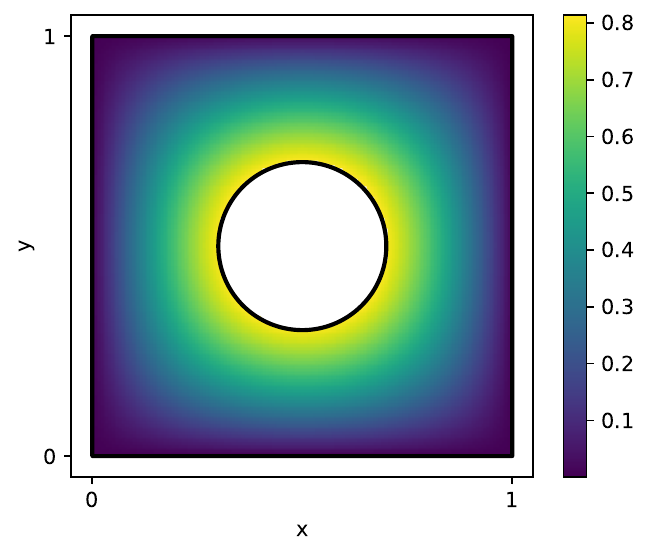}
\caption{}
\end{subfigure}\hfill
\begin{subfigure}[t]{0.4\textwidth}
\centering
\includegraphics[width=\textwidth]{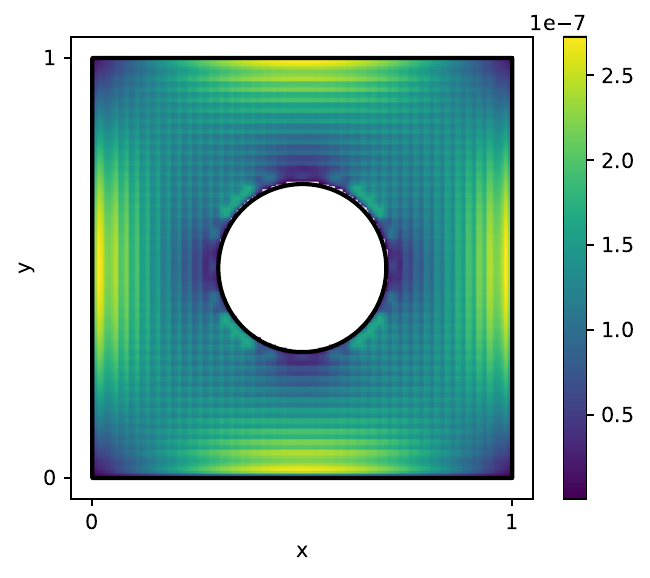}
\caption{}
\end{subfigure}
\caption{Plots of approximate solution $u_h$, (a), and spatial absolute error $\vert u-u_h\vert$, (b) for the hybrid square-hole: the outer boundary is fitted and the circular inner boundary is implicit.}
\label{fig:solution_plots_square_hole}
\end{figure}

\begin{figure} [h]
\centering
\begin{subfigure}[t]{0.4\textwidth}
\centering
\includegraphics[width=\textwidth]{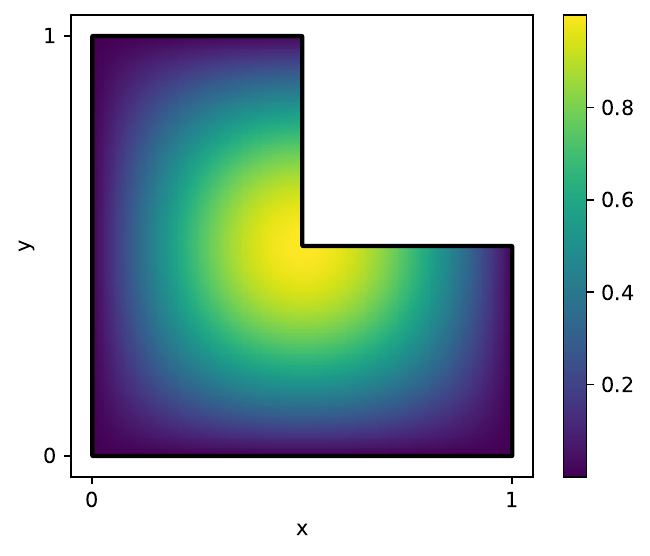}
\caption{}
\end{subfigure}\hfill
\begin{subfigure}[t]{0.4\textwidth}
\centering
\includegraphics[width=\textwidth]{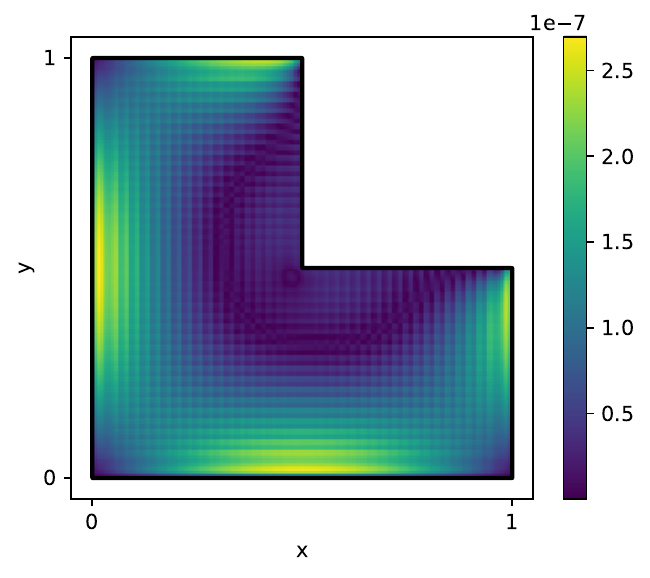}
\caption{}
\end{subfigure}
\caption{Plots of approximate solution $u_h$, (a), and spatial absolute error $\vert u-u_h\vert$, (b) for the fully fitted L-shaped domain.}
\label{fig:solution_plots_l-shape}
\end{figure}

Figure~\ref{fig:boundary_hierarchy} illustrates the locations at which the boundary conditions are imposed. Figure~\ref{fig:solution_plots_square_hole} and Figure~\ref{fig:solution_plots_l-shape} shows the computed solutions and spatial absolute errors for the hybrid and fitted L-shaped domains. In both cases, errors are evaluated on an independent grid.

The fitted L-shaped computation uses \texttt{gbmsc\_pde} \cite{el_harrak_2026_gbmsc_pde}. We retain only complete subgrids contained in the closure of the domain and impose Dirichlet data at Cartesian boundary sources. This fitted comparison uses its own retention rule, which differs from the trimming rule used for implicit domains. Since the manufactured solution is smooth, the example examines the boundary representation without testing the approximation of a corner singularity.

\FloatBarrier

\section{Conclusion}
\label{sec:conclusion}

We have developed the exterior-row replacement (ERR) framework for grid-based multinode Shepard method on implicit planar domains. The proposed idea introduces boundary completion methods that preserve the tensor product structure of retained Cartesian subgrids. The ERR strategy retains the exterior sources required by complete subgrids that intersect the boundary. It imposes the PDE collocation only at interior sources. In contrast, boundary collocation points are selected as off-source points. The full column rank is characterized by the action of the boundary rows in the interior PDE nullspace. Polynomial reproduction gives a lower bound on this nullspace through sampled harmonic polynomials, and the completion result states the conditions required when fitted rows are prescribed. A separate estimate bounds the smallest singular value using the interior action, the boundary action, and their coupling.

The numerical results support ERR-OS and ERR-QR for the tested smooth configurations. Among the 12 refinement cases, ERR-P and ERR-U are classified as full column rank in eight and five cases, respectively. ERR-OS$(1.5)$ and ERR-QR receive this classification in all twelve. Using a common candidate cloud, ERR-QR requires approximately one third fewer boundary equations than ERR-OS and gives comparable sampled value errors. However, ERR-OS is more accurate in some trimming cases, so these tests do not establish a universal preference between the two completion methods.

\section*{Code availability}
A preliminary Python implementation of the proposed method will be made available in the GitHub repository
\url{https://github.com/e-anouar/gbmsc-pde-implicit}.

\bibliographystyle{unsrt}
\bibliography{references}

\end{document}